\documentclass[12pt]{amsart}
\begin{document}

\def\tag#1{\eqno(#1)}
\def\tit#1{\par\medskip\noindent{\bf #1}\par\smallskip\noindent}
\def\ci#1{_{{}_{\scriptstyle #1}}}
\def\D{\Cal D}
\def\R{\mathbb R}
\def\C{\mathbb C}
\def\M{\widetilde M}
\def\e{\varepsilon}
\def\s{\sigma}
\def\a{\alpha}
\def\d{\delta}
\def\dim{\operatorname{dim}}
\def\f{\varphi}

\title[Cotlar inequalities for Calder\'{o}n-Zygmund operators]
{Weak type estimates and Cotlar inequalities for 
Calder\'{o}n-Zygmund operators
in nonhomogeneous spaces}

\author{F. Nazarov, S. Treil, A. Volberg}
\address{Department of Mathematics, Michigan State University,
East Lansing, Michigan, 48824}
\email[Nazarov]{fedja@math.msu.edu}
\email[Treil]{treil@math.msu.edu}
\email[Volberg]{volberg@math.msu.edu}
\curraddr[Volberg]{Mathematical Sciences Research Institute,
1000 Centennial Drive, Berkeley, CA 94707-5070}

\thanks{Partially supported by the NSF
grant DMS 9622936, binational Israeli-USA grant BSF 00030, and research programs at MSRI in the Fall of 1995 and in the Fall of 1997.}

\begin{abstract}
In the paper we consider Calder\'{o}n-Zygmund operators in
nonhomogeneous spaces. We are going to prove the analogs of classical
results for homogeneous spaces. Namely, we prove that a
Calder\'{o}n-Zygmund operator is of weak type if it is bounded in
$L^2$. We also prove several versions of Cotlar's inequality for
maximal singular operator. One version of Cotlar's inequality (a
simpler one) is proved in Euclidean setting, another one in a more
abstract setting when Besicovich covering lemma is not available. We
obtain also the weak type of maximal singular operator from these
inequalities.
\end{abstract}

\maketitle

Let $\mu$ be a measure on $\C$ satisfying the Ahlfors condition
$$
\mu(B(x,r))\le r\text{\qquad for every }x\in\C, r>0
$$
(as usual, $B(x,r):=\{y\in\C: |x-y|<r \}$).

Let $K(x,y)$ be a Calderon-Zygmund kernel, i.e.
$$
|K(x,y)|\le\frac1{|x-y|}\leqno1)
$$
and
$$
\!\! 
|K(x,y)-K(x',y)|,|K(y,x)-K(y,x')|\le\frac{|x-x'|^{\e}}
{|x-y|^{1+\e}}
\text{ whenever }|x-x'|\le \frac12|x-y|.\!\!\leqno2)
$$

\tit{Theorem 1.}
Assume that the operator
$$
T\f(x):=\int_\C K(x,y)\f(y)\,d\mu(y)
$$
(actually it would be better to write $T(\f d\mu)$ instead of
$T\f$ but we hope that the reader will forgive us for a little bit
inconsistent notation)
is bounded in $L^2(\mu)$. Then
for any (signed) measure $\nu$ on $\C$ and for any $t>0$ one has
$$
\mu\{|T\nu|>t\}\le\frac{A||\nu||}{t}.
$$
\medskip

For the Cauchy integral operator the weak type inequality was 
obtained first by X. Tolsa in [T1]. His method was specific for the 
Cauchy kernel and used the results from [Me]. Later the weak type 
estimate for general Calder\'{o}n-Zygmund operators was proved in 
[NTV2], where also the weak type estimate for maximal singular 
operator were obtained along with certain Cotlar's inequalities for 
the maximal singular operator. Another proof of the weak type estimate 
for the maximal singular operator (and yet another Cotlar inequality)
appeared in [T2].

Here and in the following we denote by $A$ a constant
depending on $||T||\ci{L^2(\mu)\to L^2(\mu)}$ and $\e$ only
(notice that we have already put the constants in the Ahlfors
condition and in the properties of the Calder\'{o}n-Zygmund
kernel $K(x,y)$ to 1: it can always be achieved by multiplication of
the measure and the kernel by a sufficiently small positive constant,
which does not change anything in the problem).

\tit{Proof.}
First of all we shall need the following standard
\tit{Obvious Lemma:}
If $\eta$ is a signed measure concentrated in some disk $B(x,\rho)$
and such that $\eta(\C)=0$, then
$$
\int_{\C\setminus B(x,2\rho)}|T\eta|d\mu\le A_1||\eta||.
$$
\tit{Proof of Obvious Lemma:}
For any $y\in\C\setminus B(x,2\rho)$ we have
\begin{equation*}
\begin{split}
|T\eta(y)|=\Bigl|\int_{B(x,\rho)}K(y,x')d\eta(x')\Bigr|=
\Bigl|\int_{B(x,\rho)}[K(y,x')-K(y,x)]d\eta(x')\Bigr|\le
\\
\le
\int_{B(x,\rho)}|K(y,x')-K(y,x)|\,d|\eta|(x')\le
||\eta||\frac{\rho^\e}{|x-y|^{1+\e}}.
\end{split}
\end{equation*}
It remains to note that
$$
\int_{\C\setminus B(x,2\rho)}\frac{\rho^\e}{|x-y|^{1+\e}}d\mu(y)\le 
1+\tfrac1\e
$$
due to the Ahlfors condition.

In order to proceed, we need to define the maximal operator
$$
T^\sharp \f(x):=\sup_{r>0}|T_r\f(x)|
$$
where
$$
T_r\f(x):=T(\f\chi\ci{\C\setminus B(x,r)})(x)=
\int_{\C\setminus B(x,r)}K(x,y)\f(y)d\mu(y).
$$
We shall also need the Hardy-Littlewood maximal function operator
$$
M\f(x):=\sup_{r>0}\frac{1}{\mu(B(x,r))}\int_{B(x,r)}|\f(y)|d\mu(y).
$$
The crux of the proof is the following
\tit{Key lemma:}
For any measurable set $F\in\C$  and any $x\in\C$
$$
T^\sharp\chi\ci F(x)\le MT\chi\ci F(x)+A_2
$$
\medskip
This lemma was inspired by (the readable part of) [DM] and
we had a strong temptation to attribute it to Guy David and Pertti 
Mattila.
Finally we have suppressed this temptation and claim it to be
our own result, but by no means we insist on the reader's
doing the same. We will postpone the proof of the key lemma
for a while and now let us derive the theorem
from it.

Without loss of generality we may assume that $\nu$ is a finite
linear combination of unit point masses with positive coefficients, 
i.e. that
$$
\nu=\sum_{j=1}^N \a_j\d_{x_j},
$$
and what is more, the $\mu$-measure of any circumference centered
at any of the points $x_j$ is equal to $0$: an arbitrary measure
$\nu$ can be obtained as a week limit of such measures; the
discussion of why one always may pass to the limit is equivalent
to the discussion of how to define $T\nu$ rigorously in the general
case. We leave this headache to the reader to take care of.

Also we can always assume that $||\nu||=\sum_j\a_j=1$ (this is just 
the
matter of normalization) and that the kernel
$K(x,y)$ is real-valued (otherwise consider $Re\,K$ and $Im\,K$ 
separately).
Thus we have to prove that
$$
\mu\{|T\nu|>t\}\le \frac At.
$$

If $\mu(\C)<\frac1t$, there is nothing to do. Otherwise let
$B(x_1,r_1)$ be a disk such that 
$\mu(E_1):=\mu(B(x_1,r_1)=\frac{\a_1}{t}$,
$B(x_2,r_2)$ be a disk such that $\mu(E_2):=\mu(B(x_2,r_2)\setminus
B(x_1,r_1)=\frac{\a_2}{t}$, and so on:
$B(x_j,r_j)$ is a disk such that
$$
\mu(E_j):=
\mu\Bigl(B(x_j,r_j)\setminus
\bigcup_{\ell=1}^{j-1}B(x_\ell,r_\ell)
\Bigr)=\frac{\a_j}{t},
$$
(the definition signs define $E_j$).
Let $E=\bigcup_j E_j=\bigcup_j B(x_j,r_j)$
Clearly,
$$
\mu(E)=\frac1t.
$$
Now let us compare $T\nu$ to $t\sum_j\chi\ci{\C\setminus B(x_j,2r_j)}
\cdot T\chi\ci{E_j}=:t\sigma$ outside $E$.
We have
$$
T\nu-t\s=\sum_j\f_j
$$
where
$$
\f_j=\a_j T\d_{x_j}- t\chi\ci{\C\setminus B(x_j,2r_j)}\cdot 
T\chi\ci{E_j}.
$$
Note now that
$$
\int_{\C\setminus E}|\f_j|d\mu\le \int_{\C\setminus 
B(x_j,2r_j)}|T[\a_j
\d_{x_j}-t\chi\ci{E_j}d\mu]|d\mu + \int_{B(x_j,2r_j)\setminus 
B(x_j,r_j)}
\a_j|T\d_{x_j}|d\mu.
$$
But, according to Obvious Lemma, the first integral does not exceed
$$
A_1||\a_j
\d_{x_j}-t\chi\ci{E_j}d\mu||=2A_1 \a_j,
$$
while $|T\d_{x_j}|\le r_j^{-1} $ outside $B(x_j,r_j)$ and therefore
the second integral is not greater than
$a_j r_j^{-1}\mu(B(x_j,2r_j))\le 2\a_j$.
Finally we conclude that
$$
\int_{\C\setminus E}|T\nu-t\s|d\mu\le 2(A_1+1)\sum_j\a_j=2(A_1+1)
$$
and thereby the absolute value of the difference $T\nu-t\s$ does not
exceed $t$ everywhere on $\C\setminus E$, except, maybe, a set of
measure $\frac{2(A_1+1)}{t}$. To accomplish the proof of the theorem,
it is enough to show that for sufficiently large $A_3$
$$
\mu\{|\s|>A_3\}\le \frac2t.
$$
We will do it by the standard Stein-Weiss duality trick.
Assume that the inverse inequality holds. Then either
$
\mu\{\s>A_3\}> \frac1t,
$
or
$
\mu\{\s<-A_3\}> \frac1t.
$
Assume for definiteness that the first case takes place
and choose some set $F\subset\C$ of measure exactly $\frac1t$
such that $\s>A_3$ everywhere on $F$.
Then, clearly,
$$
\int_\C \s\chi\ci F d\mu>\frac{A_3}{t}.
$$
On the other hand this integral can be computed as
$$
\sum_j \int_\C [T\chi\ci{E_j}]\cdot\chi\ci{F\setminus 
B(x_j,2r_j)}\,d\mu
=\sum_j \int_\C \chi\ci{E_j}\cdot [T^*\chi\ci{F\setminus 
B(x_j,2r_j)}]\,d\mu
$$
Note that for every point $x\in E_j\subset B(x_j,r_j)$
\begin{equation*}
\begin{split}
|T^*\chi\ci{F\setminus B(x_j,2r_j)}(x)-T^*\chi\ci{F\setminus 
B(x,r_j)}(x)|
=\Bigl|\int_{B(x_j,2r_j)\setminus B(x,r_j)}K(y,x)d\mu(y)  \Bigr|
\le
\\
\le r_j^{-1}\mu(B(x_j,2r_j)\le 2
\end{split}
\end{equation*}
and thereby
$$
|T^*\chi\ci{F\setminus B(x_j,2r_j)}(x)|\le (T^*)^\sharp \chi\ci 
F(x)+2\le
MT^*\chi\ci F(x)+A_2+2
$$
according to the key lemma.

Hence
$$
\int_\C \s\chi\ci F d\mu\le (A_2+2)\mu(E)+\int_\C\chi\ci E\cdot 
MT^*\chi\ci F
d\mu.
$$
But the first term equals $\frac{A_2+2}{t}$ while the second one does 
not
exceed
$$
\|\chi\ci E\|\ci{L^2(\mu)}
\|MT^*\chi\ci F\|\ci{L^2(\mu)}\le \frac{1}{t}||M||\ci{L^2(\mu)\to 
L^2(\mu)}
||T^*||\ci{L^2(\mu)\to L^2(\mu)}.
$$
Recalling that
$$
||M||\ci{L^2(\mu)\to L^2(\mu)}\le 100
$$
for an arbitrary measure $\mu$, we see that it is enough to take
$$
A_3=A_2+2+100||T||\ci{L^2(\mu)\to L^2(\mu)}
$$
to get a contradiction.
The theorem is proved and it remains to prove the key lemma.

\tit{Proof of Key Lemma.}
Let $x\in\C, r>0$. Consider the sequence of disks $B_j:=B(x,3^jr)$
and the corresponding sequence of measures $\mu_j:=\mu(B_j)$
($j=0,1,\dots)$. Note
that it is impossible that $\mu_j>9\mu_{j-1}$ for every $j\ge 1$
(this would imply that $\mu(B(x,R))>\text{const}\,R^2$ for large $R$,
which contradicts the Alfors condition). Therefore there exists
the smallest positive integer $k$ for which $\mu_k\le 9\mu_{k-1}$.
Put $R=3^{k-1}r$. Note that
\begin{equation*}
\begin{split}
|T_r\chi\ci F(x)-T\ci{3R}\chi\ci F(x)|\le
\int_{B(x,3R)\setminus B(x,r)}\frac{d\mu(y)}{|x-y|}
\le\sum_{j=1}^k\frac{\mu_j}{3^{j-1}r}\le
\\
\le\sum_{j=1}^k\frac{\mu_k 9^{j+1-k}}{3^{j-1}r}=
27\cdot 9^{-k}\mu_k r^{-1}\sum_{j=1}^k3^{j}\le
81\frac{\mu_k}{3R}\le 81.
\end{split}
\end{equation*}
And that's basically all, because now it is enough to pick
up any standard proof based on the doubling condition to get the 
desired
estimate for $T\ci{3R}\chi\ci F(x)$ (recall that $\mu(B(x,3R)\le 
9\mu(B(x,R))\,!!!)$

One of such standard ways is to compare
$T\ci{3R}\chi\ci F(x)$ to the average
$$
U\ci R(x):=\frac{1}{\mu(B(x,R))}\int_{B(x,R)}T\chi\ci Fd\mu
$$
(the quantity, which is clearly bounded by $MT\chi\ci F(x)$).

We have
\begin{equation*}
\begin{split}
T\ci{3R}\chi\ci F(x)-U\ci{R}(x)=
\\
=
\int_{F\setminus B(x,3R)} T^*[\d_x-\tfrac
{1}{\mu(B(x,R))}\chi\ci{B(x,R)}d\mu]\,d\mu-
\frac{1}{\mu(B(x,R))}\int_\C\chi\ci{B(x,R)}\cdot T\chi\ci{F\cap
B(x,3R)}\,d\mu.
\end{split}
\end{equation*}
The first term does not exceed $2A_1$ according to Obvious Lemma,
while the second can be estimated by
\begin{equation*}
\begin{split}
\frac{1}{\mu(B(x,R))}     
||T||\ci{L^2(\mu)\to L^2(\mu)}
||\chi\ci{B(x,R)}||\ci{L^2(\mu)}
\cdot ||\chi\ci{F\cap
B(x,3R)}||\ci{L^2(\mu)}\le
\\
\le
\frac{1}{\mu(B(x,R))}     
||T||\ci{L^2(\mu)\to L^2(\mu)}\sqrt{\mu(B(x,R))}\sqrt{\mu(B(x,3R))}
\le 3||T||\ci{L^2(\mu)\to L^2(\mu)}.
\end{split}
\end{equation*}

Note at last, that we can restrict ourselves to "good" disks
$B(x,r)$ in the definition of $M\f(x)$, namely, to the disks 
satisfying
$\mu(B(x,3r))\le 81\mu(B(x,r))$ (we can replace
81 by 9 so far, but we will really need this larger constant in what 
follows).
Thus the geometry
of the space (the Besicovich maximal function theorem) is not involved
in the proof: the Vitali covering theorem (that's why we 
used the stretching factor $3$ instead of more natural $2$
in the construction) is more than enough
to show the boundedness of the restricted maximal function operator
(we will denote it by $\M$)
in $L^2(\mu)$.

As usual, by interpolation we conclude that $T$ is bounded in 
$L^p(\mu)$
for every $1<p\le2$ and then, by duality, this result automatically
extends to all $p\in(1,+\infty)$.

Now we are ready to prove the boundedness of the maximal operator 
$T^\sharp$
in all spaces $L^p(\mu)$ with $1<p<+\infty$.

Let us introduce one more maximal function
$$
M'\f(x):=\sup_{r>0} r^{-1}\int_{B(x,r)}|\f|d\mu.
$$
Note that $M'$ is a bounded operator in $L^p(\mu)$ for all $p>1$
(provided that $\mu$ satisfies the Alfors condition, of course) and 
that
again the Vitali covering theorem is enough to prove this.

Let now $\beta>0$. Define
$$
M'\ci\beta\f:= \bigl\{M'[|\f|^\beta]\bigr\}^{1/\beta};\qquad
\M\ci\beta\f:= \bigl\{\M[|\f|^\beta]\bigr\}^{1/\beta}
$$

Clearly, both $M'\ci\beta$ and $\M\ci\beta$ are bounded operators
 in $L^p(\mu)$ for all $p>\beta$ 

The boundedness of the maximal operator
$T^\sharp$ follows now from 

\tit{Theorem 2.}
For any $\beta>1$
$$
T^\sharp\f\le \M T\f + B \M\ci\beta\f + B' M'\ci\beta\f
$$
where the constants $B$ and $B'$ depend on $\e$, $\beta$ and
$||T||\ci{L^2(\mu)\to L^2(\mu)}$ only.

\tit{Proof.}
It is just a minor modification of the proof of Key Lemma.
Let again $r>0$. Define $B_j$ and $\mu_j$ as before,
but let now $k$ be the smallest positive integer
for which $\mu_{k+1}\le 81\mu_{k-1}$ (i.e. we look now two
steps forward when checking for the doubling).
Let $R=3^{k-1}r$ exactly as before.
We have  
$$
|T_r\f(x)-T\ci{3R}\f(x)|\le
\int_{B(x,3R)\setminus B(x,r)}\frac{|\f(y)|d\mu(y)}{|x-y|}
\le\sum_{j=1}^k\frac{I_j}{3^{j-1}r}
$$
Where 
$$
I_j:=\int_{B_j}|\f|d\mu\le 
\Bigl\{\int_{B_j}|\f|^\beta\Bigr\}^{1/\beta}\mu_j^\gamma\le
M'\ci\beta\f(x)\,(3^jr)^{1/\beta}\mu_j^\gamma
$$
where $\gamma=1-\frac1\beta$.
Now we may continue the estimate:
$$
\sum_{j=1}^k\frac{I_j}{3^{j-1}r}
\le 3M'\ci\beta\f(x)
\sum_{j=1}^k\Bigl[\frac{\mu_j}{3^{j-1}r}\Bigl]^\gamma
\le 3M'\ci\beta\f(x)
\Bigl[\frac{\mu_k}{3^{k-1}r}\Bigl]^\gamma
\sum_{j=1}^k 3^{(j+2-k)\gamma}\le B'_1 M'\ci\beta\f(x).
$$
So, again, we need only to estimate $T\ci{3R}\f(x)$.
As before, define  
$$
U\ci R(x):=\frac{1}{\mu(B(x,R))}\int_{B(x,R)}T\f d\mu
$$
(the quantity, which is clearly bounded by $\M T\f(x)$)
and write
\begin{equation*}
\begin{split}
T\ci{3R}\f(x)-U\ci{R}(x)=
\\
=
\int_{\C\setminus B(x,3R)} T^*[\d_x-\tfrac
{1}{\mu(B(x,R))}\chi\ci{B(x,R)}d\mu]\,\f\,d\mu-
\frac{1}{\mu(B(x,R))}\int_\C\chi\ci{B(x,R)}\cdot T[\f\chi\ci{
B(x,3R)}]\,d\mu.
\end{split}
\end{equation*}
We leave it as an excercise for the reader to show
that under the conditions of Obvious Lemma
$$
\int_{\C\setminus B(x,2\rho)}|T\eta|\,|\f|\,d\mu\le
A_4 M'\f(x)\,||\eta||,
$$
and thereby the first term does not exceed
$2A_4 M'\f(x)\le 2A_4 M'\ci\beta\f(x)$.
As to the second term, we know that $T$ acts in $L^\beta(\mu)$
and therefore this term is bounded by
$$
\frac{1}{\mu(B(x,R))}||T||\ci{L^\beta(\mu)\to L^\beta(\mu)}
||\chi\ci{B(x,R)}||\ci{L^{\beta'}(\mu)}\cdot ||\f\chi\ci{
B(x,3R)}||\ci{L^\beta(\mu)}
$$
where ${\beta'}:=\frac{\beta}{\beta-1}$ is the conjugate exponent to 
$\beta$.
But clearly
$$
||\chi\ci{B(x,R)}||\ci{L^{\beta'}(\mu)}=\bigl\{\mu(B(x,R)
\bigr\}^{1/{\beta'}}
$$
while
$$
||\f\chi\ci{
B(x,3R)}||\ci{L^\beta(\mu)}\le
\M\ci{\beta}\f(x) \bigl\{\mu(B(x,3R)
\bigr\}^{1/\beta} \le
\M\ci{\beta}\f(x) \bigl\{81\mu(B(x,R)
\bigr\}^{1/\beta}                  ,
$$
and, finally, the second term is not greater than
$$
81||T||\ci{L^\beta(\mu)\to L^\beta(\mu)}\M\ci{\beta}\f(x)
$$
(we can write $\M$ instead of $M$ because due to the choice of
$k$ we have $\mu(B(x,9R))\le 81 \mu(B(x,R))$ and thereby both
balls $B(x,R)$ and $B(x,3R)$ are used in the definition of $\M$ at
the point $x$).

Now, to complete the "classical $L^p$-theory", it remains to prove
that the maximal operator $T^\sharp$ satisfies a weak type 
$L^1$-estimate,
namely, that

\tit{Theorem 3.}
For any (signed) measure $\nu$ on $\C$ and for any $t>0$ one has
$$
\mu\{|T^\sharp\nu|>t\}\le\frac{A||\nu||}{t}.
$$
\medskip

\tit{Proof.}
Let us again assume that $\nu=\sum_{j=1}^N\a_j\d_{x_j}$, $\a_j>0$,
$||\nu||=\sum_j\a_j=1$.
Fix some $t>0$ and carry out the construction of Theorem 1.
Define
$$
m(x):=\sum_j\frac{\a_j}{r_j}\chi\ci{B(x_j,10r_j)}.
$$
Let again
$$
\f_j:=T[\a_j\d_{x_j}]-t\chi\ci{\C\setminus B(x_j,2r_j)}\cdot 
T\chi\ci{E_j}.
$$
Let now $x\in\C\setminus E$, $r>0$.
Put
$$
\s_r(x):=\sum_{j}\chi\ci{\C\setminus B(x_j,2r_j)}\cdot 
T\chi\ci{E_j\setminus
B(x,r)}.
$$
Our first aim is to estimate the difference $T_r\nu(x)-t\s_r(x)$.
Clearly, it can be represented as
$$
\sum_{j=1}^N\bigl\{
T_r[\a_j\d_{x_j}](x)-t
\chi\ci{\C\setminus B(x_j,2r_j)}(x)\cdot T\chi\ci{E_j\setminus
B(x,r)}(x)=:\sum_{j=1}^N\psi_j(x)
\bigr\}.
$$
There are $3$ possibilities for each $j$:

\leftline{1)\quad $B(x_j,r_j)\cap B(x,r)=\emptyset$}

Then $\psi_j(x)=\f_j(x)$ and thereby the sum of such terms
does not exceed $$\sum_{j=1}^N|\f_j(x)|.$$

\leftline{2)\quad $B(x_j,r_j)\cap B(x,r)\not=\emptyset$, $r\le r_j$}

Since at any rate $|\psi_j(x)|\le \frac{2\a_j}{r_j}$ and since
for every such $j$ we have $x\in B(x_j,2r_j)\subset B(x_j, 10r_j)$,
we conclude that the sum of all such terms is not greater
than $2m(x)$.

\leftline{2)\quad $B(x_j,r_j)\cap B(x,r)\not=\emptyset$, $r> 2r_j$}

In this case we shall use the estimate
$|\psi_j(x)|\le \frac{2\a_j}{r}=\frac{2t\mu(E_j)}{r}$.
Since for every such $j$ one has $E_j\subset B(x_j,r_j)\subset 
B(x,3r)$
and since $E_j$ are pairwise disjoint, we obtain that the sum of such 
terms
can be estimated from above by
$$
\frac{2t\mu(B(x,3r)}{r}\le 6t.
$$
Now let us recall that (due to Obvious Lemma)
$$
\int_{\C\setminus E}\sum_{j=1}^N|\f_j|d\mu\le 2A_1
$$
and note that
$$
\int_\C md\mu\le 10.
$$
Thus, it remains only to check that for some large $A_4$
$$
\mu\{x\in\C: \s^\sharp(x)>A_4\}\le \frac{A_4}{t},
$$
where
$$
\s^\sharp(x):=\sup_{r>0}|\s_r(x)|.
$$
The proof will resemble a lot the proof of Key Lemma.
We will show that
$$
\s^\sharp(x)\le \M\s(x)+A_5 + \frac{2}{t}m(x)\tag*
$$
The last  term is due to the cutting factors $\chi\ci{\C\setminus
B(x,2r_j}$ in the definition of $\s$. Without them there would be no
difference from Key Lemma at all either in the result or in the proof.

Let us demonstrate first of all that this inequality implies the 
needed
estimate for $\s^\sharp$. In order to do that it is enough to
check that
$$
\int_\C|\s|^2d\mu\le \frac At.
$$
We will do it by duality again. Let $\f\in L^2(\mu)$. We have
$$
\int_{\C}\s\,\f\,d\mu=\sum_{j=1}^N \int_{E_j} T^*[\chi\ci{\C\setminus
B(x_j,2r_j)}\f] d\mu =:\sum_j I_j.
$$
But for every $z\in E_j$
$$
T^*[\chi\ci{\C\setminus
B(x_j,2r_j)}\f](z)= (T^*)_{r_j}[\chi\ci{\C\setminus
B(x_j,2r_j)}\f](z)-\int_{B(x_j,2r_j)\setminus B(z,r_j) 
}K(y,z)\f(y)d\mu(y).
$$
The absolute value of the first term does not exceed
$(T^*)^\sharp\f(z)$ while the second one is not greater
than $2M'\f(z)$. Hence (recall that both $(T^*)^\sharp$ and $M'$ are
bounded in $L^2(\mu)$)
$$
\Bigl|
\sum_j I_j
\Bigr|    \le
\int_E[(T^*)^\sharp\f+2M'\f]d\mu\le \sqrt{\mu(E)}\cdot 
A||\f||\ci{L^2(\mu)}\le
\frac{A||\f||\ci{L^2(\mu)}}{\sqrt t},
$$
and we are done.

Clearly, if $J$ is some subset of the index set $\{1,\dots,N\}$
and
$$
\s^{(J)}:=
\sum_{j\in J}\chi\ci{\C\setminus B(x_j,2r_j)}\cdot T\chi\ci{E_j},
$$
the same reasoning yields the estimate
$$
\int_{\C}|\s^{(J)}|^2d\mu\le A\sum_{j\in J}\mu(E_j).
$$

The first step in the proof of $(*)$ is the same as in Key Lemma:
instead of $\s_r(x)$ it is enough to consider $\s_{3R}(x)$ with
some $R$ satisfying $\mu(B(x,9R))\le 81\mu(x,R))$. Here the factors
$\chi\ci{\C\setminus B(x_j,2r_j)}$ can only help and we leave this 
step
to the reader as an excercise.

Now consider the difference
$$
\s_{3R}(x)-\frac{1}{\mu(B(x,R))}\int_{B(x,R)}\s \,d\mu.
$$
It can be represented as the sum
$$
\sum_{j=1}^N\bigl\{
\chi\ci{\C\setminus B(x_j,2r_j)}(x)\cdot T\chi\ci{E_j\setminus 
B(x,3R))}(x)-
\frac{1}{\mu(B(x,R))}\int_{B(x,R)}\chi\ci{\C\setminus B(x,2r_j)}\cdot
T\chi\ci{E_j}=:\sum_j D_j,
d\mu
\bigr\}.
$$
We will split the index set $\{1,\dots,N\}$ into several subsets:

\leftline{1)\qquad $J_1:=\{j:B(x,3R)\cap B(x_j, 2r_j)=\emptyset\}$}

For every $j\in J_1$ we have
$$
D_j=\int_{E_j} T^*[\d_{x}-\frac{1}{B(x,R)}\chi\ci{B(x,R)}d\mu]d\mu
$$
and, since $E_j$ are pairwise disjoint and since $E_j\subset 
\C\setminus
B(x,3R)$ for every $j\in J_1$, Obvious Lemma shows that
$$
\Bigl|\sum_{j\in J_1}D_j\Bigr|\le 2A_1.
$$

\leftline{2)\qquad $J_2:=\{j:B(x,3R)\cap B(x_j, 2r_j)\ne\emptyset,\
R<2r_j\}.$} For any $j\in J_2$ we have $x\in B(x_j, 10r_j)$ and, due 
to
the trivial estimate $D_j\le 2\frac{\mu(E_j)}{r_j}$, we
conclude that
$$
\Bigl|\sum_{j\in J_2}D_j\Bigr|\le 2\frac{m(x)}{t}.
$$

\leftline{3)\qquad $J_3:=\{j:B(x,3R)\cap B(x_j, 2r_j)\ne\emptyset
,\ R\ge 2r_j\}$}

Note that for every $j\in J_3$ one has $E_j\subset B(x,9R)$.
Now observe that
$$
\bigl|
\chi\ci{\C\setminus B(x_j,2r_j)}(x)\cdot T\chi\ci{E_j\setminus 
B(X,3R))}(x)
\bigr|\le \frac{\mu(E_j)}{3R}
$$
and thereby
$$
\sum_{j\in J_3}
\bigl|
\chi\ci{\C\setminus B(x_j,2r_j)}(x)\cdot T\chi\ci{E_j\setminus 
B(X,3R))}(x)
\bigr|\le \sum_{j\in J_3}\frac{\mu(E_j)}{3R}\le
\frac{\mu(B(x,9R))}{3R}\le 3.
$$
On the other hand,
\begin{multline}
\Bigl|\frac{1}{\mu(B(x,R)}\int^{B(x,R)}\s^{(J_3)}d\mu\Bigr|\le
\frac{1}{\mu(B(x,R)}
\sqrt{\mu(B(x,R)}
\Bigl\{\int^{B(x,R)}\bigl|\s^{(J_3)}\bigl|^2d\mu\Bigr\}^{1/2}\le
\\
\le
\sqrt{\frac{A\sum_{j\in J_3}\mu(E_j)}{\mu(B(x,R)}}\le
\sqrt{\frac{A\mu(B(x,9R))}{\mu(B(x,R)}}\le 9\sqrt A.
\end{multline}
Thus
$$
\Bigl|\sum_{j\in J_3}D_j\Bigr|\le 3+9\sqrt A
$$
and we are through.

Let $M$ stands for a usual maximal function, that is let
$$
M\f(x) := \sup_{r>0}\mu(B(x,r))^{-1} \int_{B(x,r)}|\f| d\mu.
$$

This function is larger than $M^{'}\phi$ introduced earlier. But in 
cases when this maximal function has weak type $(1,1)$ or is bounded 
in some other sense (which will be the case for Euclidean space 
setting by Besicovich covering lemma) we may give an easier proof of 
the weak boundedness of $T^{\sharp}$. The proof of the theorem itself 
does not require any geometry or any covering properties. It is the 
application that is sensitive to the type of covering theorems. 
Theorem below was obtained also by Xavier Tolsa in [T2]. He uses  
the earlier version [NTV2] to do that. Actually he uses a certain idea from [DM] but in
the interpretation given in [NTV2].

In [T2] this result is  used give a full
characterization of positive measures
$\mu$  on the  complex plane such that the truncations of Cauchy integral of any 
complex measure would pointwisely converge with respect to $\mu$. In 
particular, if $\mu$ is such that the Cauchy integral operator is 
bounded in $L^2(\mu)$ then this pointwise convergence takes place. 
These remarkable result finishes the series of works [NTV1], [T1], 
[NTV2], [T2], in which the nonhomogeneous Calder\'{o}n-Zygmund theory 
was considered. 

Coming back to the weak type estimates for maximal singular
operator we first prove the following Cotlar's inequalities.
 
\tit{Theorem 4.}
$$
T^\sharp\f\le A_p \{M |T\f|^p\}^{1/p} + B_p M\f,
$$
if $p \in (0,1)$. Also
$$
T^\sharp\f\le A M T\f + B M\f 
$$
where the constants $A$ and $B$ depend on 
$||T||\ci{L^2(\mu)\to L^2(\mu)}$ only.
\tit{Proof.}
It is just a minor modification of the proof of Key Lemma.
Let again $r>0$.
First we are proving the following form of Cotlar's inequality. Fix 
$p < 1$.

Clearly the next inequality implies both inequalities in Theorem 4.
$$
|(T_r F)(y)|^{p} \le M(|TF|^{p})(y)+C (MF(y))^{p}.
$$
Fix $r>0$ and let $R$ be the number from Key Lemma. 
Let us prove first that 
$$
|(T_r F)(y)| \le |(TF)(x)| + |(TF\chi_{B(y, 3R)})(x)|+ C MF(y)
$$
for any $ x \in B(y,R)$.

Suppose this is proved. Then we rais it to the power $p$ and average 
over $B(y,R)$ with respect to measure $\mu$. We will need the 
following estimate
$$
\int_{B(y,R)}|(TF\chi_{B(y, 3R)})(x)|^p d\mu(x) \le C 
\mu(B(y,R))^{1-p}(\int_{B(y,3R)}|F| d\mu)^p
$$
Dividing by $\mu(B(y,R))$, and noticing that the choice of $p$ gives
$\mu(B(y,3R)) \le 9 \mu(B(y,R))$ we get that the averaging of 
$|(TF\chi_{B(y, 3R)})(x)|^p$ over $B(y,R)$ is bounded by 
$$
9^p C( \mu (B(y, 3R))^{-1}\int_{B(y,3R))}|F| d \mu)^p \le 9^p C 
(MF(y))^p
$$. 
Then we get our first inequality, and we are done.

So we are left to prove two things

1) $|(T_r F)(y)| \le |(TF)(x)| + |(TF\chi_{B(y, 3R)})(x)|+ C MF(y)$,

2) $\int_{B(y,R)}|(TF\chi_{B(y, 3R)})(x)|^p d\mu(x) \le C 
\mu(B(y,R))^{1-p}(\int_{B(y,3R)}|F|d\mu)^p$.

As for the first one, its proof repeats Key Lemma.  Define $B_j$ and 
$\mu_j$ as in the proof of Key Lemma. Then
$$
|(T_r F)(y) -T_{3R}F)(y)| \le \sum_{j=1}^{k} 
(3^{-j}R)^{-1}\mu(B(y,3R))3^{-2j}
\mu(B_j)^{-1} \int_{B_j}|F| d\mu \le C MF(y).
$$
On the other hand, for $x \in B(y,R)$ obviously
$$
|T_{3R}F)(y)-T_{3R}F)(x)|| \le C MF(y)
$$
and $|T_{3R}F)(x)| \le |(TF)(x)| + |T(F \chi_{B(y,3R)})(x)|$. Thus 
the first inequality is proved.

Now the second one follows from the following elementary lemma.

\tit{Sublemma }
Let $\f \in L^{1, \infty}$ and let $p \in (0,1)$. For a measurable 
set $B$ we have
$$
\int_{B}|\f|^p \le C_p (\mu(B))^{1-p} (\| \f \|_{L^{1, \infty}})^p.
$$
\tit{Proof} It is an obvious use of the formula
$$\int_{B}|\f|^p d \mu = \int_{0}^{\infty} p s^{p-1}\mu\{y \in B: \f 
> s\} ds
\le \int_{0}^{\infty} p s^{p-1}\min(\mu(B), 
s^{-1}\|\f\|_{L^{1,\infty}}) ds
$$
Now we apply Sublemma to $\f =T(F\chi_{B(y, 3R)})$, and we get 2) 
immediately if we use already proven fact that $T$ maps $L^1$ to 
$L^{1, \infty}$.

Theorem 4 is completely proved now.

{\bf Remark}. It is clear from the proof that the maximal function in front of $|T\f|$
(or $|T\f|^p$ term is a restricted maximal function $\tilde{M}$. Meaning that it is
taken only over the discs with doubling property. One need not Besicovich covering lemma
to give the weak type estimate for such restricted maximal function. But the second
maximal function--the one in front of $\f$--is the usual unrestricted centered maximal
function. This is why Theorem 4 is applicable only in Euclidean spaces.

Let us show why Theorem 4 gives the weak type estimate for $T^{\sharp}$. We use the
first inequality of Theorem 4. The second term is estimated as usual, and by Besicovich
covering lemma it is in $L^{1,\infty}$. To estimate the first term let us notice that
$$
\mu\{x: MF(x) > t\} \le Ct^{-1}\int_{F >t/2}F d\mu
$$
for positive functions $F$. To notice that we just split $F$ to $F \chi_{\{x: F\le
t/2\}}$ and $F \chi_{\{x: F > t/2\}}$. Then apply the usual weak type estimate to $M(F
\chi_{\{x: F > t/2\}})$. After this remark let us apply it to $F:= |T\f|^p$ and $t:=
2 s^p$. Using Sublemma again (and also using twice the fact that $T$ maps $L^1$ to
$L^{1,\infty}$) we get
$$
\mu\{x: M(|T\f|^p) > 2 s^p \}\le Cs^{-p}\int_{|T\f| > s}|T\f|^p d\mu
$$
$$
\le C s^{-p} \mu\{ |T\f| >s\}^{1-p}\|\f\|_1^p \le C s^{-p}s^{p-1}\|\f\|_1^p
\|\f\|_1^{1-p},
$$
and we are done.

\tit{References}

[DM] G. David, P. Mattila, Removable sets for Lipschitz harmonic 
functions in the plane, Prepublication, ORSAY, 1997, No. 31, pp. 1-59.

[Me] M. Melnikov, Analytic capacity: discrete approach and curvature of the measure.
Mat. Sbornik, {\bf 186} (1995), No. 6, pp. 827-846.

[NTV1] F. Nazarov, S. Treil, A. Volberg, Cauchy integral and 
Calder\'{o}n-Zygmund operators on nonhomogeneous spaces. IMRN 
International Math. Res. Notes, 1997, {\bf No. 15}, pp. 703-726.

[NTV2] F. Nazarov, S. Treil, A. Volberg, Calder\'{o}n-Zygmund 
operators in nonhomogeneous spaces. Preprint, June 1997, pp. 1-10.

[T1] X. Tolsa, $L^2$-boundedness of the Cauchy integral operator for 
continuous measures. Preprint, 1997, pp. 1-24.

[T2] X. Tolsa, Cotlar's inequality and existence of principal values 
for the Cauchy integral without doubling conditions. Preprint, 1997, 
pp. 1-31.

\end{document}